# A Heterogeneous Multiscale Method for Power System Simulation Considering Electromagnetic Transients

Kaiyang Huang[1], *Student Member, IEEE*, Min Xiong[1], *Student Member, IEEE,* Yang Liu[2], *Member*, *IEEE*, Kai Sun[1], *Senior Member, IEEE,*
Feng Qiu[2], *Senior Member, IEEE*

[1]University of Tennessee, Knoxville, TN, USA
[2]Argonne National Laboratory, Lemont, IL, USA
khuang12@vols.utk.edu, mxiong3@vols.utk.edu, yliu161@anl.gov, kaisun@utk.edu, fqiu@anl.gov

*Abstract*—Traditional dynamic security assessment faces challenges as power systems are experiencing a transformation to inverter-based-resource (IBR) dominated systems, for which electromagnetic transient (EMT) dynamics have to be considered. However, EMT simulation is time-consuming especially for a large power grid because the mathematical model based on detailed component modeling is highly stiff and needs to be integrated at tiny time steps due to numerical stability. This paper proposes a heterogeneous multiscale method (HMM) to address the simulation of a power system considering EMT dynamics as a multiscale problem. The method aims to accurately simulate the macroscopic dynamics of the system even when EMT dynamics are dominating. By force estimation using a kernel function, the proposed method automatically generates a macro model on the fly of simulation based on the micro model of EMT dynamics. It can flexibly switch between the micro- and macro-models to capture important EMT dynamics during some time intervals while skipping over other time intervals of less interest to achieve a superior simulation speed. The method is illustrated by a case study on a two-machine EMT model to demonstrate its potential for power system simulation.

## I. Introduction

In recent decades, time-domain simulation has been widely applied to dynamic security assessment and transient stability analysis for power systems [1]. However, with increasing Inverter-Based-Resources (IBRs) dominating the system, power system dynamics become more complex to include both electromechanical dynamics and electromagnetic transient (EMT) dynamics, and thus positive-sequence phasor-based models are not adequate to provide accurate results. Indeed, several subsynchronous oscillations observed in recent years can hardly be captured by phasor models [2]. EMT simulation tools such as the Electromagnetic Transients Program (EMTP) can accurately simulate a three-phase power network model even with distributed transmission line parameters, and become more popular in analyses of fast dynamics for the desired accuracy. However, EMT simulations are extremely slow because the instantaneous voltage and current waveforms around the fundamental 60-Hz need to be computed at a step size of milliseconds or shorter for numerical stability even if the system is approaching a steady-state condition. Several straightforward approaches were proposed to speed up EMT analysis by a transformation of the system model and responses into frequency domain such as the Hilbert Transform and Fast Fourier Transform [3],[4]. Such methods have limitations such as the aliasing effect with fast dynamics to accumulate their errors, and it is inconvenient to analyze the error bound. Thus, two questions arise for a power system having IBRs: How to define fast and slow dynamics of different timescales? How to efficiently simulate a multi-timescale power system having both fast and slow dynamics?

There have been works on co-simulations of both EMT and phasor models by interfacing two simulation tools[5-7]. However, very few papers have concerned the methodology for multiscale simulations of power systems considering EMT dynamics. This paper aims to fill the gap and develops a Heterogenous Multiscale Method (HMM) for simulations of power systems having two timescales. The HMM [8-10] provides a general methodology to efficiently simulate a stiff system, which relies on the coupling of different time scales, i.e., a micro-model can supply the necessary information for the macro-model. The method is based on time-averaging techniques from the perturbation theory [11] for analyzing stiff dynamical systems. The proposed novel HMM solver could flexibly and efficiently exchange information and switch between micro- and macro-models of different timescales with errors bounded. The main advantage of the proposed approach is that the effective force is estimated through sophisticatedly defined smooth kernels. By using a convolution operator with kernels, the macro-model can capture the dynamics of the power system correctly.

The rest of this paper is organized as follows. Section II derives and studies the formulae on the HMM for a stiff system. Then, the HMM solver is introduced, and a kernel function is implemented to capture the effective force. Section III introduce the EMT model of a two-machine system under 0DQ global reference frame. Section IV presents case studies on the EMT model of a two-machine system. Conclusions are summarized in section V.

## II. Heterogeneous Multiscale Method

This section explains the motivation and basic idea of the

This work was supported by the Advanced Grid Modeling (AGM) program of U.S. DOE Office of Electricity under grant DE-OE0000875.

proposed method using a stiff dynamical system. The influence of a smooth kernel is studied on such a system. Then, the HMM solver is proposed and the algorithm for multiscale power system simulation is built.

*A. Introduction of a stiff system*

A stiff ordinary differential system can be defined by:

$$\frac{dx_\varepsilon}{dt} = f_\varepsilon(x_\varepsilon, t) \quad (1)$$

where $x_\varepsilon(t): \mathbb{R}^+ \mapsto \mathbb{R}^d$ is the state variable, $d$ is the dimension of the system. Subscript $\varepsilon$ indicates that the system has a faster timescale characterized by a small positive number $\varepsilon$. From [10], the eigenvalues of the system should satisfy, $\exists k_0 \geq 1$ such that:

$$\begin{cases} \Re \lambda_\varepsilon(j) \leq C_1, & 1 \leq j \leq d \\ |\lambda_\varepsilon(j)| \leq C_2, & 1 \leq j \leq k_0 \leq d \\ C_3 \leq \varepsilon |\lambda_\varepsilon(j)| \leq C_4, & k_0 < j \leq d \\ \min_{j_1, j_2} |\lambda_\varepsilon(j_1) - \lambda_\varepsilon(j_2)| > \rho > 0 \end{cases} \quad (2)$$

where from $C_1$ to $C_4$ are all constants and $\lambda_\varepsilon(j)$ denotes the $j^{th}$ eigenvalues of the Jacobian of the vector field in (1), also $j_1 \leq k_0$ and $j_2 > k_0$. Typically, a stiff system is illustrated roughly by eigenvalues in different scales, (2) can explain such a property precisely.

Consider an example with two dimensions $(x_1, x_2)$:

$$\begin{cases} \varepsilon \dot{x}_1 = f_1(x_1, x_2, t) \\ \dot{x}_2 = f_2(x_1, x_2, t) \end{cases} \quad (3)$$

where $f_1$ and $f_2$ are all smooth functions, $x_2$ is a slow variable compared with $x_1$ when $\varepsilon > 0$ is small enough. Such a stiff system can be approximated by Differential Algebraic Equations (DAE) by applying perturbation theory. Indeed, when $\varepsilon$ converges to 0, the system becomes:

$$\begin{cases} 0 = f_1(x_1, x_2, t) \\ \dot{x}_2 = f_2(x_1, x_2, t) \end{cases}. \quad (4)$$

Hence the fast dynamics in the first equation converge to a slow manifold. In the power system area, the DAE model for transient stability comes from the convergence of power flow equations.

In addressing issues for multiscale dynamics in power system simulations, now suppose that there exists an effective macro model behind the fast dynamics:

$$\frac{d}{dt}X = \bar{f}(X, t) \quad (5)$$

where $X(t): \mathbb{R}^+ \mapsto \mathbb{R}^q$ represents the slow variables derived from (1) as $\varepsilon$ converges to 0. Here we emphasize the dimension of slow variables to be different from the micro system (1). Notice that different from the simple example (4), the macro model is generated automatically on the fly of simulation and hence may not have a closed form, so we only need it to exist. For power system simulation, we are interested in some slow dynamics such as the envelope of fast-oscillating state variables, the energy of the system, or relative angles, they may usually have such forms. Then, based on (5), the effective force at $t$ can be defined by [10]:

$$\bar{f}(t) = \lim_{\delta \to 0} \left[ \lim_{\varepsilon \to 0} \frac{1}{\delta} \int_t^{t+\delta} f_\varepsilon(\tau) d\tau \right]. \quad (6)$$

It is almost impossible to solve (6) explicitly for a nontrivial problem, so a strategy is to use time averaging techniques to evaluate the effective force from a micro model, i.e. original system (1). Indeed, by selecting a proper kernel, a convolution operation with kernel function can achieve this goal. Similar to [10], we consider a kernel function $K(t) \in \aleph^{p,q}(I)$ which means: $K(t) \in C_c^q(\mathbb{R})$ with supp$(K) = I$ and

$$\int_\mathbb{R} K(t) t^r dt = \begin{cases} 1, & r = 0 \\ 0, & 1 \leq r \leq p \end{cases}. \quad (7)$$

For normalization, $K_\eta(t)$ is defined to denote the kernel function after scaling, i.e.,

$$K_\eta(t) \doteq \frac{1}{\eta} K(\frac{t}{\eta}).$$

Based on the above definitions, the estimation of an effective force by the convolution using a kernel function $K_\eta^{p,q}$ with $\eta = \eta(\varepsilon) \to 0$ as $\varepsilon \to 0$ converges[8], i.e.

$$\tilde{f} = K_\eta^{p,q} * f_\varepsilon = K_\eta^{p,q} * (\bar{f} + g_\varepsilon) \to \bar{f} \text{ as } \varepsilon \to 0. \quad (8)$$

where $g_\varepsilon(t)$ represents fast dynamics that need to be averaged in the macro model (5). For detailed proof, please check **Theorem 2.3** and **2.7** for the dissipative and oscillatory cases [10].

**Remark 1.** For a stiff system (1), we define a macro model (5) of interest, the effective force (6) of (5) usually cannot have a closed form. However, by time averaging using a kernel function, (6) can be obtained from the micro-model.

*B. Algorithm of the heterogeneous multiscale method*

Based on the knowledge introduced in the previous section, a frame can be built systematically to address this EMT model. Suppose the simulation is started at $t_n$, the simulation algorithm is presented below:

*Step 1* Estimation of macro effective force:
  a) Reconstruct information from the macro-model by
  $$x_0 = RX^{(n)}$$
  b) Solve micro-model (EMT) based on the micro solver:
  $$\frac{dx_\varepsilon}{dt} = f_\varepsilon(x_\varepsilon, t) \text{ with } x_\varepsilon(t_n) = x_0 \text{ for } t \in [t_n, t_n + \eta]$$
  c) Apply time averaging to the micro-model (EMT):
  i) Using (8) to estimate the local macro force, i.e.
  $$\bar{f}(t_n + \Delta t) \approx \tilde{f}_n = \tilde{f}(t_n + \Delta t) = K_\eta^{p,q} * f_\varepsilon(t_n + \Delta t)$$
  ii) Obtain macro information from the micro-model (EMT)
  $$X^* = Q(x_\varepsilon(t_n + \Delta t))$$

*Step 2* Evolve the macro dynamics $X^{n+1}$ for the next step $t_{n+1}$:

$$X^{(n+1)} = \sum_{k=m}^{n} A_k X^{(k)} + H \sum_{k=m}^{n} B_k \tilde{f}_n + CX^*$$

*Step 3* Let $x_0 = RX^{(n+1)}$, then repeat the whole process.

where $R$ and $Q$ are reconstruction and compression operators transforming between micro- and macro-state variables, e.g., $Q = \|\cdot\|_{L2}$ represents the general energy form of micro-state variables. $\eta$ denotes the simulation interval of the micro-model at $t_n$, $X^{(n)}$ denotes slow variables at step $n$. Also, $\Delta t$ determines the force evaluation position as $t_n + \Delta t$. In ii) in step 1, all vector fields which have the form like $f(x(t_n+ \Delta t))$ are written in short as $f(t_n+ \Delta t)$. In step 2), a general form of a multistep method is considered to illustrate the process of the macro simulation, where $A_k$, $B_k$, and $C$ in $\mathbb{R}^{q \times q}$ represent coefficient matrices and $H$ is the step size for macro simulation. This algorithm provides a general frame to capture the dynamics we are interested in power system simulation, i.e., the evaluation of macro dynamics is achieved on the fly of simulation for the micro-model. In this paper, we consider $R$, and $Q$ as identity operators, which means the slow dynamics of state variables in (1) are what we prefer to observe during the simulation. By repeating this process, the effective force (6) can be predicted from the micro model (1) instead of our definition such as the positive-sequence phasor-based model which may cause errors in some cases as pointed out in [2]. The basic logic is that most fast dynamics $f_\varepsilon$ are not cores of our study, by running the micro-model (1) with a time interval $\eta$, enough information is obtained to support the simulation for the macro-model (5) by a convolution with a proper kernel function in (8), which depends on the problems we address.

**Remark 2**. Indeed, to assess the overall stability and dynamic performance of the system, we do not focus on state variables of fast dynamics since they either decay quickly right after a switch or enter steady-state periods such as instantaneous voltage and current waveforms around the fundamental 60-Hz. Some methods shift such dynamics to a low frequency [3][4].

*C. The influence of a kernel function in HMM*

In the previous section, a robust and adaptive HMM algorithm is introduced, and the role of a kernel function is analyzed in this section. By local averaging against a kernel function, fast dynamics can be averaged but long-term dynamics can be preserved.

In general, a convolution operator with a kernel is a linear continuous functional for a function. The kernel applied in this paper is unimodal and isotropic, which means the kernel has radial symmetry and should be invariant under the rotation operation. Such a kernel function could always be represented as a kind of metric, indeed, there exists a smooth function $h \in C_c^q(\mathbb{R})$ such that

$$K(t,s) = h(\|t - s\|).$$

With the abuse of notation, kernel $K$ is written as:

$$K(t,s) = K(t-s) \qquad (9)$$

as the kernel only depends on the difference. Then, the integral of a function $f(s)$ against the kernel becomes convolution naturally

$$\tilde{f}(t) = \int_\mathbb{R} K(t-s) f_\varepsilon(s) ds = (K * f_\varepsilon)(t), \qquad (10)$$

if we identify $K(t-s)$ in (9) as a single variable function, then (9) has the same form defined in (7). The convolution operation in (10) evaluates the function $f_\varepsilon$ at $t_n+\Delta t$, and $\|t_n+\Delta t - t\|$ represents the distance between two points, then the kernel can be recognized as a weighted function used to average $f_\varepsilon$ depending on the chosen metric, i.e., the feature we want to keep. The redundant component $g_\varepsilon$ can be filtered during this process.

### III. DYNAMIC MODEL UNDER 0DQ REFERENCE FRAME

In this section, we adopt a common practice to establish the dynamic model into 0DQ frame reference which can transform balanced sinusoidal dynamics into nearly constant dynamics.

*A. Model of a round rotor synchronous generator*

Consider a 6$^{th}$-order EMT model:

$$\begin{aligned}
\frac{d\delta}{dt} &= \Delta\omega_r \\
\frac{d\Delta\omega_r}{dt} &= \frac{\omega_0}{2H}(p_m - p_e - D\frac{\Delta\omega_r}{\omega_0}) \\
\frac{d\lambda_{fd}}{dt} &= e_{fd} - \frac{r_{fd}}{L_{lf}}(\lambda_{fd} - \lambda_{ad}) \\
\frac{d\lambda_{1d}}{dt} &= -\frac{r_{1d}}{L_{1dl}}(\lambda_{1d} - \lambda_{ad}) \\
\frac{d\lambda_{1q}}{dt} &= -\frac{r_{1q}}{L_{1ql}}(\lambda_{1q} - \lambda_{aq}) \\
\frac{d\lambda_{2q}}{dt} &= -\frac{r_{2q}}{L_{2ql}}(\lambda_{2q} - \lambda_{aq}) \\
\frac{di_{abc}}{dt} &= -\frac{1}{L_{abc}''}(v_{abc} - P_{ark}^{-1}v_{0dq}'' + R_s i_{abc} + \frac{dL_{abc}''}{dt}i_{abc})
\end{aligned} \qquad (11)$$

where $\delta$, $\omega_r$, $\Delta\omega_r$, $H$, $D$, $p_m$, and $p_e$ denote the rotor angle (rad), rotor angle speed (rad/s), rotor angle speed deviation(rad/s), inertial constant (s$^2$/rad), damping constant (p.u.), mechanical power (p.u.), and electrical power (p.u.), respectively; $\omega_0$ is the nominal frequency of the system, i.e., 60 Hz; $\lambda_{fd}$, $\lambda_{1d}$, $\lambda_{1q}$, and $\lambda_{2q}$ denote flux linkages (p.u.) of filed winding, damper winding at $d$-axis, the first damper winding at $q$-axis, the second damper winding at $q$-axis, respectively; $r_{fd}$, $r_{1d}$, $r_{1q}$, and $r_{2q}$ are resistances (p.u.) of those four windings; $L_{fdl}$, $L_{1dl}$, $L_{1ql}$, and $L_{2ql}$ are leakage inductances (p.u.) of those four windings; $e_{fd}$ is the field voltage (p.u.); $i_{abc}$ and $v_{abc}$ are three phase terminal currents and voltages interfacing with the grid side; $R_s$ is the corresponding stator resistance. We also consider a 1$^{st}$-order governor and exciter for each generator in the case study.

*B. Model of a three-phase RLC branch under 0DQ reference frame*

A two-machine system shown in Fig.1 is introduced to explain the model of the RLC branch. The dynamic model can be transformed into a global 0DQ reference frame as below:

$$\begin{aligned}
\frac{di_1}{dt} &= L_{T1DQ}^{-1}(v_1 - v_3) + \tilde{P}(\theta)i_1 \\
\frac{di_2}{dt} &= L_{T2DQ}^{-1}(v_2 - v_4) + \tilde{P}(\theta)i_2 \\
\frac{di_4}{dt} &= L_{2DQ}^{-1}(v_4 - R_{2DQ}i_4) + \tilde{P}(\theta)i_4 \\
\frac{di_7}{dt} &= L_{LineDQ}^{-1}(v_3 - v_4 - R_{LineDQ}i_7) + \tilde{P}(\theta)i_7 \\
\frac{dv_3}{dt} &= C_{LineDQ}^{-1}(i_1 - i_7) + \tilde{P}(\theta)v_3 \\
\frac{dv_4}{dt} &= C_{LineDQ}^{-1}(i_2 - i_4 + i_7) + \tilde{P}(\theta)v_4
\end{aligned} \qquad (12)$$

where $i_1$, $i_2$, $i_4$, and $i_7$ denote currents through G$_1$, G$_2$, load 2 and

T-line, $v_3$, $v_4$ are voltages at capacitors in T-line. All elements in the network are transformed by Park transformation:

$$L_{T1DQ} = P(\theta)L_{T1}P(\theta)^{-1} \quad L_{T2DQ} = P(\theta)L_{T2}P(\theta)^{-1}$$
$$L_{2DQ} = P(\theta)L_2 P(\theta)^{-1} \quad R_{2DQ} = P(\theta)R_2 P(\theta)^{-1}$$
$$L_{LineDQ} = P(\theta)L_{Line}P(\theta)^{-1} \quad R_{LineDQ} = P(\theta)R_{Line}P(\theta)^{-1}$$
$$C_{LineDQ} = P(\theta)C_{Line}P(\theta)^{-1}$$

And Park matrix adopted in this paper is:

$$P(\theta) = \frac{2}{3}\begin{bmatrix} 1/2 & 1/2 & 1/2 \\ \cos(\theta) & \cos(\theta - \frac{2\pi}{3}) & \cos(\theta + \frac{2\pi}{3}) \\ -\sin(\theta) & -\sin(\theta - \frac{2\pi}{3}) & -\sin(\theta + \frac{2\pi}{3}) \end{bmatrix}$$

$$\tilde{P}(\theta) = \omega_r \frac{dP}{dt}(\theta)P(\theta)^{-1}$$

In (12) all currents, voltages as well as resistors, inductors, and capacitors are under the 0DQ reference frame. Each element is diagonalized by Park matrix, $L_{T1}$ and $L_{T2}$ denote inductors of two transformers, $L_2$ and $R_2$ are used to model the load level, and $L_{Line}$, $R_{Line}$, $C_{Line}$ represent parameters of the transmission line. Notice that load 1 will be tripped in the case study, and also has the same form as load 2, so the model of load 1 isn't listed in (12).

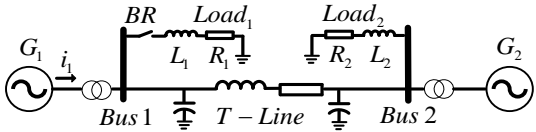

Fig. 1. A two-machine system

Note that in (12), all state variables have three components, e.g., $i_1 = [i_{10}, i_{1D}, i_{1Q}]^T$, which are in the global 0DQ reference instead of the local 0dq reference for each generator or IBR. Without loss of generality, we could choose $\theta$ and $\omega_r$ from the local 0dq frame of $G_1$ as the choice is arbitrary [12]. Then, all quantities in the network are rotating together with $G_1$. The transformation between different frames can be achieved by a rotation matrix, see details in [12].

## IV. CASE STUDIES

### A. HMM on the two-machine system

The performance of HMM is demonstrated in this section. Consider a fault that happens at $t=3$s which leads BR to open in Fig.1, i.e., load 1 is permanently tripped at $t=3$s, and HMM is adopted to simulate the dynamics. The whole process can be divided into three periods: steady state period $T_1$ before disturbance, $T_2$ the period during fault and part of post fault period in which micro dynamics dominate the system, $T_3$ the period in which long-term dynamics are important. We claim that the strategy adaptively switches between different models during the simulation. Then, during $T_3$, only long-term dynamic is preferred, to speed up the simulation, HMM is applied in this process since there are still some fast dynamics left during this period.

We set our numerical simulation for a two-machine system as follows. As discussed in the algorithm of HMM, $H$ is used to denote the macroscopic step size, $h$ denotes the step size of micro simulation, $\eta$ denotes the kernel support size, and the time window of micro simulation is $W = 2\eta$. In this paper, the Gaussian kernel is applied to handle the discrete convolution with the micro vector field $f_\varepsilon$. The Gaussian kernel has the form:

$$K(t,s) = \frac{1}{\sqrt{2\pi\sigma^2}}e^{-\frac{(t-s)^2}{2\sigma^2}}. \quad (13)$$

Identify $t-s$ as a single variable, then (13) has the same form as (7). Note that the choice of $\sigma$ depends on the micro window $W$. It would be a good strategy to take $\sigma = \eta/3$ as $K$ taking $t-s$ out of this domain is near 0. Pick the middle point as the evaluation point which can provide high accuracy due to the symmetric kernel (13), i.e., $\Delta t = \eta$.

In this case, the micro model is EMT model introduced at III. We don't define macro model explicitly. The HMM-FE-rk4 solver is applied in this case, which means RK4 is adopted during the micro simulation and Forward Euler (FE) scheme is used for macro simulation. Detailed information is shown below. In this case, set $T_2$ as 3s to 3.1s, as the fast dynamics is damped very fast, typically the high-frequency component in voltage, which is much shorter than 0.1 seconds, also we want to observe the fast dynamics for some state variables during $T_2$.

TABLE I. PARAMETERS OF HMM SOLVER

| $\sigma$ | $\Delta t$(s) | $H$(s) | $h$(μs) |
|---|---|---|---|
| 0.0044 | 0.011 | 0.012 | 5 |
| $W$(s) | $\eta$(s) | $T_1$(s) | $T_2$(s) |
| 0.022 | 0.011 | [0,3) | [3,3.1] |
| $T_3$(s) | R&Q | Micro Solver | Macro Solver |
| (3.1,8] | I* | RK4 | FE |

* $I$ is the identical operator

Then, voltages $v_5$, $v_6$, and currents $i_4$, under 0DQ frame are shown below, RK4 scheme using 5 μs time step is adopted as the benchmark for the whole simulation.

From Fig.2, we can conclude that simulation with HMM can capture the long-term dynamics accurately with a much larger macroscopic time step $H$. Furthermore, another group of currents {$i_1$, $i_2$, $i_7$} behaves stronger oscillations due to the loss of the load, and they are shown in Fig.3, Fig.4, and Fig.5. The corresponding envelope is also compared with the original EMT which aims to show the potential of the HMM. All long-term dynamics of currents can be accurately simulated with less computational burden.

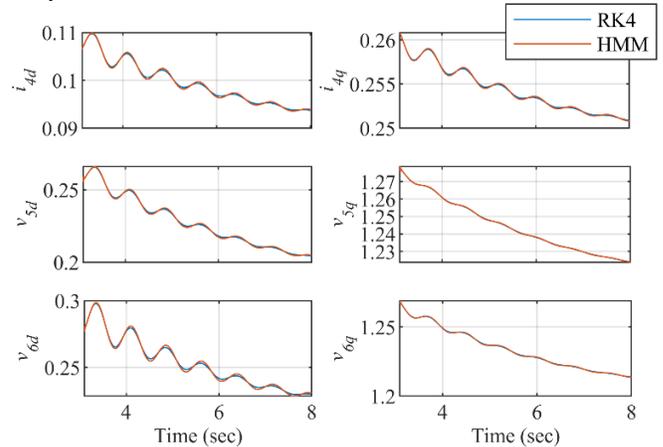

Fig. 2. Results of $i_4$, $v_5$, and $v_6$ during $T_3$.

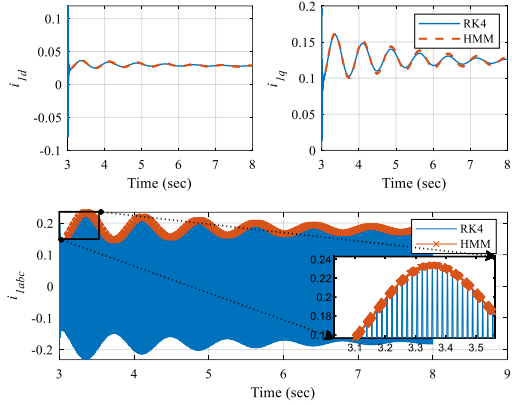

Fig. 3. Results of $i_{1d}$, $i_{1q}$, and the envelope of $i_1$ under *abc* reference frame.

From Fig. 3 it can be seen the dense part of HMM is the micro simulation which is adopted to evaluate the effective force for the macro simulation. The simulation is stable and effective for $\{i_1, i_2, i_7\}$. Also, there are several implicit versions of HMM, e.g. one can apply the implicit RK4 scheme to the micro simulation. In this case study, $T_3$ is not in a steady state as the fault is permanent. HMM can track the envelope during the transient period of the network device.

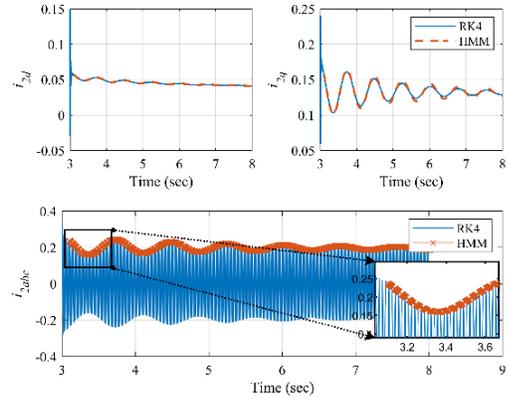

Fig. 4. Results of $i_{2d}$, $i_{2q}$, and the envelope of $i_2$ under *abc* reference frame.

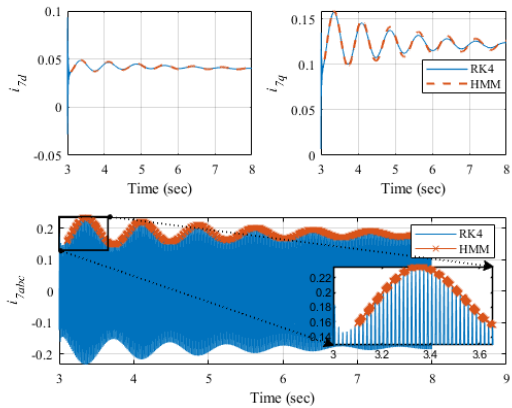

Fig. 5. Results of $i_{7d}$, $i_{7q}$, and the envelope of $i_7$ under *abc* reference frame.

The simulation is conducted in Matlab on a PC with Intel(R) Core (TM) i7-10700 CPU@2.90GHz. Regarding time performance, for such an 8-second simulation, the ground truth EMT simulation shown before is used to compare with the HMM-based simulation. HMM-based simulation is finished with 80.931170s and RK4 one is 121.932214s which implies the power of HMM, a 33.63% speedup is achieved as a much larger time step is adopted during the macro simulations.

## V. CONCLUSION AND FUTURE WORK

This is the first attempt to apply the HMM to power system simulation to achieve a balance between micro and macro dynamic in EMT simulation. The benefit of HMM demonstrated in the case study is that it provides a general adaptive flexible frame to simulate a stiff system, time performance is improved as well as accuracy is preserved. This paper opens a door for the future study of HMM in power system simulation, e.g., whether an adaptive step change scheme can be applied to increase the time performance; if a better kernel function could be found to further average the long-term dynamics.